\documentclass{article}
\usepackage{amsfonts}
\usepackage{amsmath}
\usepackage{amssymb}

\newtheorem{theorem}{Theorem}
\newtheorem{proposition}{Proposition}
\newcommand{\C}{\mathbb C}
\begin{document}
\title{A note on the location of polynomial roots}
\author{D.A. Bini\thanks{Universit\`a di Pisa, {\tt bini@dm.unipi.it}}~  and 
 F. Poloni\thanks{
Scuola Normale Superiore, Pisa, {\tt f.poloni@sns.it}}}
\maketitle
\begin{abstract}
  We review some known inclusion results for the roots of a polynomial, and
  adapt them to a conjecture recently presented by S. A. Vavasis. In
  particular, we provide strict upper and lower bounds to the distance of the
  closest root of a polynomial $p(z)$ from a given $\zeta\in\C$ such that
  $p'(\zeta)=0$.
\end{abstract}

\section{Introduction}
Recently S.A. Vavasis \cite{vavasis} has presented the following conjecture.

{\bf Conjecture} {\it There exist two universal constants $0<\iota_1\le 1\le
  \iota_2$ with the following property. Let $\xi_1,\ldots,\xi_n$ be 
 the roots of
  a degree-$n$ univariate polynomial $p(z)$. Let $\zeta_1,\ldots,\zeta_{n-1}$
  be the roots of its derivative. Define
\begin{equation}\label{eq:rho}
\rho_j=\min_{k=2,\ldots,n}\left|\frac{k!p(\zeta_j)}{p^{(k)}(\zeta_j)}
\right|^{1/k},~~~j=1,\ldots,n-1
\end{equation}
and the annuli
\[
A_j=\{z:\iota_1\rho_j\le|z-\zeta_j|\le\iota_2\rho_j\},~~~j=1,\ldots,n-1.
\]
 Then for each $i=1,\ldots,n$
\[
\xi_i\in A_1\cup \cdots\cup A_{n-1}.
\]
}

The author also refers to an unpublished communication by Giusti et Al., where
it is shown that $\iota_1$ exists and can be taken $(\sqrt 5-1)/2$ and where a
sequence of $n$-degree polynomials is given such that $\lim_n
|z-\zeta_j|/\rho_j=+\infty$ so that $\iota_2$ does not exist.

In this note we revisit some known general bounds to the roots of a polynomial
from \cite{henrici}, in particular Theorem 6.4b on pages 451,452, and Theorem
6.4e on page 454, and adapt them to the conditions of the Vavasis conjecture.
More specifically, we show that for any polynomial $p(z)$, and for any $\zeta$
such that $p'(\zeta)=0$, there exists a root $\xi$ of $p(z)$ satisfying
\[
\left|\xi-\zeta\right| \le \rho \sqrt{n/2},~~~\rho=\min_{k=2,\ldots,n}
\left(\frac{k!p(\zeta)}{p^{(k)}(\zeta)}\right)^{1/k},
\]
and that the bound is sharp since it is attained by a suitable polynomial.

We provide also some sharp lower bound to $|\xi-\zeta|$ under the condition
that $p^{(k)}(\zeta)=0$ for $k\in\Omega$, where $\Omega$ is a nonempty subset
of $\{1,2,\ldots,n-1\}$.

Moreover, we also show that $\iota_2$ does not exist by providing an example
of a sequence $\{p_n(z)\}_n$ of polynomials of degree $n+1$ having a common
root $\xi$, where the ratio $|\xi-\zeta_j^{(n)}|/\rho_j^{(n)}$ is independent
of $j$ and tends to infinity as $n^{1-\epsilon}$ for any $i=1,\ldots,n$ and
for any $0<\epsilon<1$, where $\zeta_j^{(n)}$ are the roots of $p_n'(z)$.

\section{Main results}
In this section, after providing a counterexample of the Vavasis conjecture,
we review some inclusion theorems of \cite{henrici}, which give lower
bounds and upper bounds to the distance of the roots of a polynomial
from a given complex number $\zeta$.

\subsection{Counterexample} 
Consider the monic polynomial of degree $n+1$
\[
p_n(z)=z^{n+1}-(n+1)z.
\]
Clearly $z=0$ is one of its roots, and we have
$
p_n'(z)=(n+1)(z^n-1),
$
 so that the roots $\zeta_i$ of $p_n'$ are the complex $n$-th roots of the
unity.  Define
\[
\rho^{(n,k)}=
\left\vert\frac{ k! p_n(\zeta)}{p_n^{(k)}(\zeta)}
\right\vert^{1/k},~~~~\rho^{(n)}=\min_{k=2,\ldots ,n+1} \rho^{(n,k)},
\]
where $\zeta$ stands for any $n$-th root $\zeta_i$ of 1, 
and observe that $p_n(\zeta)=-n\zeta$,
$p_n^{(2)}(\zeta)=n(n+1)\zeta^{-1}$. Therefore, for $k=2$ one has
\[
\rho^{(n)} \le \rho^{(n,2)}=\left\vert
    \frac{2!p_n(\zeta)}{p_n^{(2)}(\zeta)}\right\vert^{1/2}
    =\left|\frac{2n}{n(n+1)}\right|^{1/2}
    =\sqrt{\frac 2{n+1}}   
\]
hence $\rho_{n} \to 0$ as $n \to \infty$. Observe that this bound is
independent of the root $\zeta_i$.
The annuli $A_i$
have their centers on the unit circle and for $\iota_2$ constant,
their external radii tend to 0 as $n \to \infty$. Thus, for sufficiently large
values of $n$ they cannot contain the origin, and this contradicts the
conjecture as $z=0$ is a common root to all the polynomials $p_n(z)$.

Moreover, for $z=0$ one has 
\[
\frac{|z-\zeta_i|}{\rho^{(n)}}\ge\frac{|z-\zeta_i|}{\rho^{(n,2)}}
= (n+1)^{1/(n+1)} \sqrt{\frac{n+1}{2}}
\ge\sqrt{\frac{n+1}{2}}.
\]
 That is, the ratio $\frac{|z-\zeta_i|}{\rho^{(n)}}$ can grow as much as $\sqrt
 {n/2}$. 
For general $k$ 
one can easily get
\begin{equation}\label{controes}
\frac{|z-\zeta|}{\rho^{(n,k)}}=\left[\frac{1}{n}\binom{n+1}{k}\right]^{1/k}(n+1)^{\frac1{n+1}}\ge\left[\frac{1}{n}\binom{n+1}{k}\right]^{1/k}.
\end{equation}
Thus,  for a fixed $k$ the ratio $|z-\zeta|/\rho^{(n,k)}$ can
grow as much as $n^{1-\frac 1k}$.

\subsection{Lower bounds}
Let us recall the following result (see \cite{henrici}, Theorem 6.4b).

\begin{theorem}\label{th1}
Let $p(z)=\sum_{i=0}^na_i z^i$ be a monic polynomial of degree
$n$ and $\zeta$ any complex number. Assume $a_0\ne 0$. 
Then any root $\xi$ of $p(x)$ is such that
\begin{equation}\label{eq:lbound}
\gamma \rho
< |\xi-\zeta|,~~~\rho=\rho(\zeta)=
\min_{k=2,\ldots,n}\left| k!\frac{p(\zeta)}{p^{(k)}(\zeta)}\right|^{1/k}
\end{equation}
where $\gamma=1/2$. 
\end{theorem}

The following proof of the above theorem can be easily adjusted to the
case where $\zeta$ is a (numerical) root of some derivative of $p(z)$.

Without loss of generality we may assume $\zeta=0$. In fact, if $\zeta\ne 0$
consider $\widehat p(z)=p(z-\zeta)$ so that $\widehat p'(z)=p'(z-\zeta)$ and
$\widehat \rho(0)=\rho(\zeta)$, and reduce the case to $\zeta=0$.
 
From the definition of $\rho$ one has
\begin{equation}\label{eq:1}
\rho^k\le k!\left| \frac{p(0)}{p^{(k)}(0)}\right|=\left|\frac{a_0}{a_k}\right|.
\end{equation}
Then taking the moduli in both
sides of the equation $-a_0=a_1 \xi+a_2 \xi^2+\ldots+a_n \xi^n$ yields
\[
1\le \sum_{i=1}^n |\frac{a_i}{a_0} \xi^i|
\]
which, in view of \eqref{eq:1} provides the bound
\[
1\le \sum_{i=1}^n t^i,~~~ t=\frac {|\xi|}{\rho} ,
\]
whence
\[
1\le \frac{t-t^{n+1}}{1-t}.
\]
If  $t<1$ then we have  $1-t\le t-t^{n+1}<t$ which implies $t>1/2$.
This proves the bound $|\xi|>\frac 12 \rho$ for any root $\xi$ of
$p(z)$.

Observe that the bound is strict since the polynomial $p_n(z)=\sum_{i=1}^n
z^i-1$ has a root in the interval $(1/2,1/2(1+1/n))$ for $n\ge 2$.

The proof of Theorem \ref{th1} can be adjusted to the case where $\zeta$
satisfies some additional condition. We have the following result:

\begin{proposition}\label{th2} Assume that $\zeta$ satisfies the following
  condition
\[
  \theta^i\left|\frac{p^{(i)}(\zeta)}{i!p(\zeta)}\right|
  \le\epsilon,~~~i\in\Omega=\{i_1,\ldots,i_h\}\subset\{1,2,\ldots,n-1\}
\]
where $0\le \epsilon<1/h$,  $1\le h<n$ and $\theta$ is an
upper bound to $|\zeta-\xi_i|$ for $i=1,\ldots,n$. Then
\eqref{eq:lbound} holds where $\gamma$ is the only solution in
$(1/2,1)$ of
the equation
\begin{equation}\label{eq:equat}
(t-1)\sum_{i\in\Omega} t^i+2t-1+(1-t)h\epsilon=0.
\end{equation}
\end{proposition}
{\bf Proof.}
By following the same arguments of the proof of Theorem \ref{th1} 
with $\zeta=0$ one
obtains
\[1\le \sum_{i=1}^n\left|\frac{a_i}{a_0}\xi^i\right|\le
\sum_{i=1,n;~i\not\in\Omega}
\left|\frac{a_i}{a_0}\xi^i\right|+h\epsilon\le
\sum_{i=1,n;~i\not\in\Omega}t^i+h\epsilon.
\]
If $t<1$, replacing $\sum_{i=1,n;~i\not\in\Omega}t^i=(t-t^{n+1})/(1-t) 
-\sum_{i\in\Omega} t^i$
in the latter inequality yields
$1-t\le t-t^{n+1}-(1-t)\sum_{i\in\Omega}t^i+(1-t)h\epsilon\le 
t+(t-1)\sum_{i\in\Omega}t^i+(1-t)h\epsilon$. 
Whence, $t>\gamma$ where $\gamma$ is the only solution of 
\eqref{eq:equat} in $(1/2,1)$.
\hfill $\Box$\\[2ex]

Let us look at some specific instances of the above result. For $\epsilon=0$
the condition of the proposition turns into $p^{(i)}(\zeta)=0$ for
$i\in\Omega$. If in addition $\Omega=\{1\}$ one finds the condition
$p'(\zeta)=0$ of the Vavasis conjecture and \eqref{eq:equat} turns into
$t^2+t-1=0$ that implies $\gamma=(\sqrt 5-1)/2=0.618...$.  Weaker bounds are
obtained assuming $\epsilon=0$ and $\Omega=\{k\}$ for some $k>1$ since the
only root of the polynomial $t^{k+1}-t^k+2t-1$ in $(1/2,1)$ is lower than
$(\sqrt 5-1)/2$.

Better bounds are obtained if $\zeta$ is a root of multiplicity $h$ of
$p'(z)$; in fact, $\gamma$ is the only positive root of the polynomial
$t^{h+1}+t-1$. In particular, if $h=2$ then $\gamma=0.682... $, if $h=3$,
$\gamma=0.724...$.

If $\zeta$ is close to a root of $p'(z)$, so that the condition
$\theta|p'(\zeta)/p(\zeta)|<\epsilon$ for some ``small'' $\epsilon$ is
satisfied, then $\gamma=(\sqrt 5-1)/2 -\epsilon (1+3/\sqrt 5)+O(\epsilon^2)$.

For $\epsilon=0$ the bound in the above proposition is strict since it is
asymptotically attained by the polynomial $t^n- (t-1)\sum_{i\in\Omega}
t^i-2t+1$.  The advantage of this bound is that it allows to compute sharper
values for $\gamma$ just by solving a low degree equation if $\Omega$ is made
up by small integers.

Slightly better lower bounds can be obtained from the following known result
of \cite{henrici} which requires to compute a positive root of a polynomial of
degree $n$.

\begin{theorem}\label{th3}
Any root $\xi$ of $p(z)$ is such that $|\xi|\ge \sigma$, where
$\sigma$ is the only positive solution to the equation
$|a_0|=\sum_{i=1}^n t^i |a_i|$.
\end{theorem}

\subsection{Upper bounds}
Throughout this section we denote
\[
\rho^{(k)}=\left(
  k!p(\zeta)/p^{(k)}(\zeta)\right)^{1/k},~~~\rho=\min_k \rho^{(k)} 
\]
for a given $\zeta\in\C$.
Concerning upper bounds to the distance of a root from $\zeta$ we
recall the following result of \cite{henrici} (Theorem 6.4e, page 454).
\begin{theorem}\label{th4}
For
any $\zeta\in\C$ there exists a root $\xi$ of $p(z)$ such that
\begin{equation}\label{eq:ubound1}
|\xi-\zeta|\le \rho^{(k)}
{n\choose k}^{1/k} ,~~~k=1,\ldots,n.
\end{equation}
\end {theorem}

Observe that, for $k=2$ one has 
\begin{equation}\label{eq:ubound2}
 |\xi-\zeta|\le \rho^{(2)}\sqrt{n(n-1)/2},
\end{equation}
while
\begin{equation}\label{eq:hen}
|\xi-\zeta|\le\min_k {n\choose k}^{1/k} \rho^{(k)}\le
\max_k{n\choose
 k}^{1/k}\rho \le n\rho.
\end{equation}
The bound \eqref{eq:hen} is sharp since it is attained by the polynomial
$p(z)=(z-n)^n$ with $\zeta=0$.  In fact, it holds $\rho=\rho^{(1)}=1$ and
$p(z)$ has roots of modulus $n$.

Under the condition $p'(\zeta)=0$ the bounds \eqref{eq:ubound1},
\eqref{eq:ubound2} and \eqref{eq:hen} can be substantially improved.  In fact
we may prove the following result

\begin{proposition}\label{th:ub}
For any $\zeta\in\C$ such that $p'(\zeta)=0$ there exists a root $\xi$ of
$p(z)$ such that
\begin{equation}\label{eq:newbound1}
|\xi-\zeta|\le\left\{
\begin{array}{ll}
\rho^{(2)}\sqrt{n/2}\\[2ex]
\rho^{(3)}\sqrt[3]{n/3}\\[2ex]
\rho^{(k)} \sqrt n \left(\frac 1k \prod_{i=2}^{\lfloor k/2\rfloor}
   ( \frac 1n+\frac 1{2i-1}+\frac 1{2i-2})\right)^{1/k}&
\hbox {for}~~4\le k\le n \\[2ex]
\end{array}\right.
\end{equation}
Moreover, 
\begin{equation}\label{eq:newbound2}
|\xi-\zeta|\le \rho \sqrt \frac n2
\end{equation}
\end{proposition}

{\bf Proof.} Without loss of generality we may assume $\zeta=0$  and $a_0=1$
so that the polynomial can be written as
$p(z)=1+a_2z^2+\ldots+a_n z^n$. Recall
the Newton identities \cite{henrici}, page 455:
\[
k a_k=-s_k-\sum_{i=1}^{k-1} a_i s_{k-i}, ~~k=1,2,\ldots,
\]
where $s_k=\sum_{i=1}^n\xi_i^{-k}$ are the power sums of the reciprocal of the
roots $\xi$ of $p(z)$.
Clearly, $a_1=s_1=0$ so that for $k\ge 4$ the Newton identities turn into
\begin{equation}\label{eq:ni}
ka_k=-s_k-\sum_{i=2}^{k-2} a_i s_{k-i}, ~~k=4,5,\ldots.
\end{equation}
Let $\Delta=\min_i|\xi_i|$ so that $|s_k|\le n\Delta^{-k}$. It holds
$|2a_2|=|s_2|\le n\Delta^{-2}$, $|3a_3|=|s_3|\le n\Delta^{-3}$ and
\[
k|a_k|\le \Delta^{-k}n(1+\sum_{i=2}^{k-2} |a_i|\Delta^i),~~k\ge 4.
\]
Denoting $\gamma_k=n(1+\sum_{i=2}^{k-2} |a_i|\Delta^{i})$, for $k\ge 4$ and
$\gamma_2=\gamma_3=n$, by using the induction argument one easily finds that
\begin{equation}\label{eq:rec}\begin{array}{l}
k|a_k|\le\Delta^{-k}\gamma_k\\
\gamma_k\le\gamma_{k-1}+\frac n{k-2}\gamma_{k-2},~~k\ge 4\\
\gamma_2=\gamma_3=n.
\end{array}
\end{equation}
The above expression provides the bound
\begin{equation}\label{eq:boundp}
\Delta\le\rho^{(k)}\left(\frac{\gamma_k}k\right)^{1/k}
\end{equation}
so that it remains to give upper bounds to $\gamma_k$. 
Since $\gamma_2=\gamma_3=n$,  
from \eqref{eq:boundp} we deduce \eqref{eq:newbound1}
for $k=2,3$. For the general case $k\ge 4$, 
we express
the recurrence \eqref{eq:rec} in matrix form as
\[
\left[\begin{array}{c}\gamma_{k+1}\\ \gamma_{k}\end{array}\right]
  \le\left[\begin{array}{cc}
1&\frac n{k-1}\\1&0\end{array}\right]
\left[\begin{array}{c}\gamma_k\\ \gamma_{k-1}\end{array}\right],
\] 
where the inequality holds component-wise.  Applying twice the above
bound  yields 
\begin{equation}\label{eq:mf}
\left[\begin{array}{c}\gamma_{k+1}\\ \gamma_{k}\end{array}\right]
  \le\left[\begin{array}{cc}
1+\frac n{k-1}&\frac n{k-2}\\1&\frac n{k-2}\end{array}\right]
\left[\begin{array}{c}\gamma_{k-1}\\ \gamma_{k-2}\end{array}\right].
\end{equation}
Whence, since $\gamma_2=\gamma_3=n$, one finds that $\gamma_{2i}$ and
$\gamma_{2i+1}$ are polynomials in $n$ of degree $i$. 
Denoting
\begin{equation}\label{eq:delta}
\gamma_{2i}=n^i\delta_{2i},~~~~\gamma_{2i+1}=n^i\delta_{2i+1},
\end{equation}
we may give
upper bounds to  $\delta_k$. In fact, from 
\eqref{eq:mf} with $k=2i$ it holds
\begin{equation}\label{eq:mf1}
\left[\begin{array}{c}\delta_{k+1}\\ \delta_{k}\end{array}\right]
  \le\left[\begin{array}{cc}
\frac 1n+\frac 1{k-1}&\frac 1{k-2}\\\frac 1n&\frac 1{k-2}\end{array}\right]
\left[\begin{array}{c}\delta_{k-1}\\ \delta_{k-2}\end{array}\right].
\end{equation}
Let us denote $W_k$ the matrix in the right-hand side of \eqref{eq:mf1}, 
so that for $n\ge 4$ we have
\begin{equation}\label{eq:x}
\left[\begin{array}{c}\delta_{2i+1}\\\delta_{2i}\end{array}\right]=W_{2i}W_{2(i-1)}\cdots W_4\left[\begin{array}{c}\delta_{3}\\\delta_{2}\end{array}\right].
\end{equation}
Since for $n\ge 4$ 
  we have
$||W_k||_\infty=\frac 1n+\frac 1{k-1}+\frac 1{k-2}$, taking norms in 
\eqref{eq:x} yields 
\[
||(\delta_{2i+1},\delta_{2i})||_\infty\le\prod_{j=2}^i||W_{2j}||_\infty 
||(\delta_3,\delta_2)||_\infty\le
\prod_{j=2}^i\left(\frac1n+\frac 1{2j-1}+\frac 1{2j-2}\right),
\] 
since $||(\delta_3,\delta_2)||_\infty=||(1,1)||_\infty=1$.
In  view of \eqref{eq:boundp} and \eqref{eq:delta} 
this proves \eqref{eq:newbound1}.

In order to prove the bound \eqref{eq:newbound2}, from \eqref{eq:boundp} it is
sufficient to prove that
\begin{equation}\label{eq:aux}
\gamma_k\le k (\sqrt{\frac n 2})^k.
\end{equation} 
We prove the latter bound by induction on $k$ for $2\le k\le n$. 
For $k=2,3$, the inequality
\eqref{eq:aux} is true since $\gamma_2=\gamma_3=n$.  
Moreover, from \eqref{eq:rec} one has $\gamma_4\le\gamma_3+\frac n2
\gamma_2=n(n+2)/2$ so that \eqref{eq:aux} is satisfied also for $k=4$.
Now we assume that the
bound \eqref{eq:aux} is true for $k$ and $k-1$, where $k\ge 4$ and we prove it
for $k+1\le n$, i.e., $\gamma_{k+1}\le(k+1) (\sqrt{n/2})^{k+1}$.  From
\eqref{eq:rec} and from the inductive assumption one has $\gamma_{k+1}
=\left(\sqrt{\frac n2}\right)^{k+1}\left(k\sqrt{\frac 2n}+2\right) 
$.
Therefore it is sufficient to prove that
$k\sqrt{\frac 2n}+2\le k+1$, that is, $\sqrt{\frac n2}\ge \frac k{k-1}$
which is satisfied for $n\ge k\ge 4$.
 This completes the proof.\hfill $\Box$
\\[1ex]

Observe that the bound of Theorem \ref{th:ub} is sharp since it is attained by
the polynomial $p(z)=(z^2-m)^m$ with $\zeta=0$, where $n=2m$.  In fact,
$p'(0)=0$, $\rho=\rho^{(2)}=1$ and the roots of $p(z)$ have moduli
$\sqrt{n/2}$.

If $\zeta$ is such that $p^{(j)}(\zeta)=0$, $j=1,\ldots,h$, then from the
Newton identities one finds that $s_i=a_i=0$, $i=1,\ldots,h$ so that equation
\eqref{eq:ni} turns into
\[
ka_k=-s_k-\sum_{i=h+1}^{k-h-1} a_i s_{k-i}, ~~k\ge 2(h+1).
\]
By following the same argument used in the proof of Proposition \ref{th:ub} we 
can prove that there exists a root $\xi$ of $p(z)$ such that
\[
|\xi-\zeta|\le
\rho^{(h+i)}\sqrt[h+i]{\frac n {h+1}},~~i=1,\ldots,h+1.
\]


\begin{thebibliography}{99}
\bibitem{henrici} P. Henrici,
{\em Applied and Computational Complex Analysis},
Vol. 1, Wiley, 1974.
\bibitem{vavasis} S. A. Vavasis, A conjecture that the roots of a univariate
  polynomial lie in a union of annuli (Interim Revised Version),
  arXiv:math.CV/0606194 v3, 28 Jul 2006.
\end{thebibliography}
\end{document}